\documentclass[11pt]{article}
\usepackage[top=2.5cm, bottom=2.5cm, left=2.5cm, right=2.5cm]{geometry}
\usepackage{hyperref,color,graphicx,cite,array,booktabs,longtable}
\usepackage{amsmath,amssymb,amsfonts,amsthm}

\newtheorem{theorem}{Theorem}
\newtheorem{lemma}[theorem]{Lemma}
\newtheorem{corollary}[theorem]{Corollary}

\theoremstyle{definition}
\newtheorem{definition}{Definition}
\newtheorem{example}{Example}

\renewcommand{\-}{\!-\!}
\newcommand{\+}{\!+\!}

\begin{document}

\title{On Hosoya's dormants and sprouts%
\footnote{This work was supported and funded by Kuwait University Research Grant No. SM03/20. 
\par Corresponding author: Salem Al-Yakoob.
\par Email addresses: smalyakoob@gmail.com (Salem Al-Yakoob), akanso@hotmail.com, (Ali Kanso), 
			    dragance106@yahoo.com (Dragan Stevanovi\'c).}}
\author{Salem Al-Yakoob$^a$, Ali Kanso$^a$, Dragan Stevanovi\'c$^b$}
\maketitle
\begin{center}
$^a${\em Department of Mathematics, Faculty of Science, Kuwait University, Safat 13060, Kuwait;}
$^b${\em Mathematical Institute, Serbian Academy of Sciences and Arts, Kneza Mihaila 36, 11000 Belgrade, Serbia}
\end{center}

\medskip
\begin{abstract}
The study of cospectral graphs is one of the traditional topics of spectral graph theory.
Initial expectation by theoretical chemists G\"unthard and Primas in 1956
that molecular graphs are characterized by the multiset of eigenvalues of the adjacency matrix
was quickly refuted by the existence of numerous examples of cospectral graphs
in both chemical and mathematical literature.
This work was further motivated by Fisher in 1966 in the influential study that investigated 
whether one can ``hear'' the shape of a (discrete) drum, 
which has led over the years to the construction of many cospectral graphs.  
These findings culminated in setting the ground for the Godsil-McKay local switching and the Schwenk's use of coalescences, 
both of which were used to show (around the 1980s) that almost all trees have cospectral mates.
Recently, enumerations of cospectral graphs with up to 12 vertices by Haemers and Spence and by Brouwer and Spence 
have led to the conjecture that, on the contrary, ``almost all graphs are likely to be determined by their spectrum''. 
This conjecture paved the way for myriad of results showing that various special types of graphs are determined by their spectra. 

On the other hand, in a recent series of papers,
Hosoya drew the attention to a particular aspect of constructing cospectral graphs by using coalescences:
that cospectral graphs can be constructed by attaching multiple copies of a rooted graph 
in different ways to subsets of vertices of an underlying graph.
The principal focus of this research effort is to address the expectations and questions raised in Hosoya’s papers.
We give an explicit formula for the characteristic polynomial of such multiple coalescences,
from which we obtain a necessary and sufficient condition for cospectrality of these coalescences.
We enumerate such pairs of cospectral multiple coalescences for a few families of underlying graphs,
and show the infinitude of cospectral multiple coalescences having paths as underlying graphs,
which were deemed rare by Hosoya.


\end{abstract}

\smallskip
\begin{center}
{\em\small Dedicated to Haruo Hosoya for his 85th birthday.\\[6pt]\phantom{.}}
\end{center}

\section{Introduction}


An initial belief by chemists G\"unthard and Primas in 1956~\cite{gupr}
that the multiset of eigenvalues of adjacency matrix may characterize graphs
was quickly nullified by constructing pairs of graphs sharing the same spectrum of eigenvalues
(the so-called {\em cospectral} graphs).
The first such pair was identified by Collatz and Sinogowitz in their seminal paper~\cite{cosi}.
Furthermore, examples followed soon in both mathematical and theoretical chemistry literature.
At that time, the H\"uckel molecular orbital theory postulated that 
the energy levels of $\pi$-electrons in molecules are 
determined by adjacency eigenvalues of their molecular graphs.
Hence, the differences in the physicochemical properties of molecules that share the same spectrum
would point out to phenomena that this theory could not explain.
The same problem of the existence of cospectral graphs has also emerged 
in concert with research directions related to mathematical physics.
Kac in~\cite{kac} modeled the drum's shape in a continuous fashion and
showed that its sound is characterized by the eigenvalues of the eigenvalue problem
defined on the region of the drum membrane and its boundary.
Moreover, Fischer~\cite{fisch} modeled the shape of the drum in a discrete manner by a graph,
and posed the famous question: `Can one hear the shape of a drum?',
which was later equivalently translated to
whether a graph can be characterized by the multiset of its eigenvalues
(and Fischer then found additional examples of appropriate cospectral graphs).

\begin{figure}[ht]
\begin{center}
\includegraphics[scale=0.7]{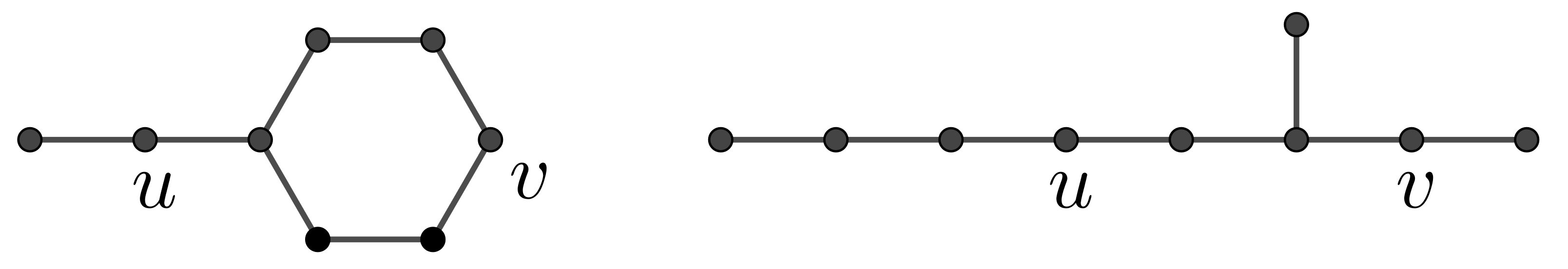}
\end{center}
\caption{Examples of graphs with isospectral vertices.}
\label{fig-1}
\end{figure}

\begin{figure}[ht]
\begin{center}
\includegraphics[scale=0.6]{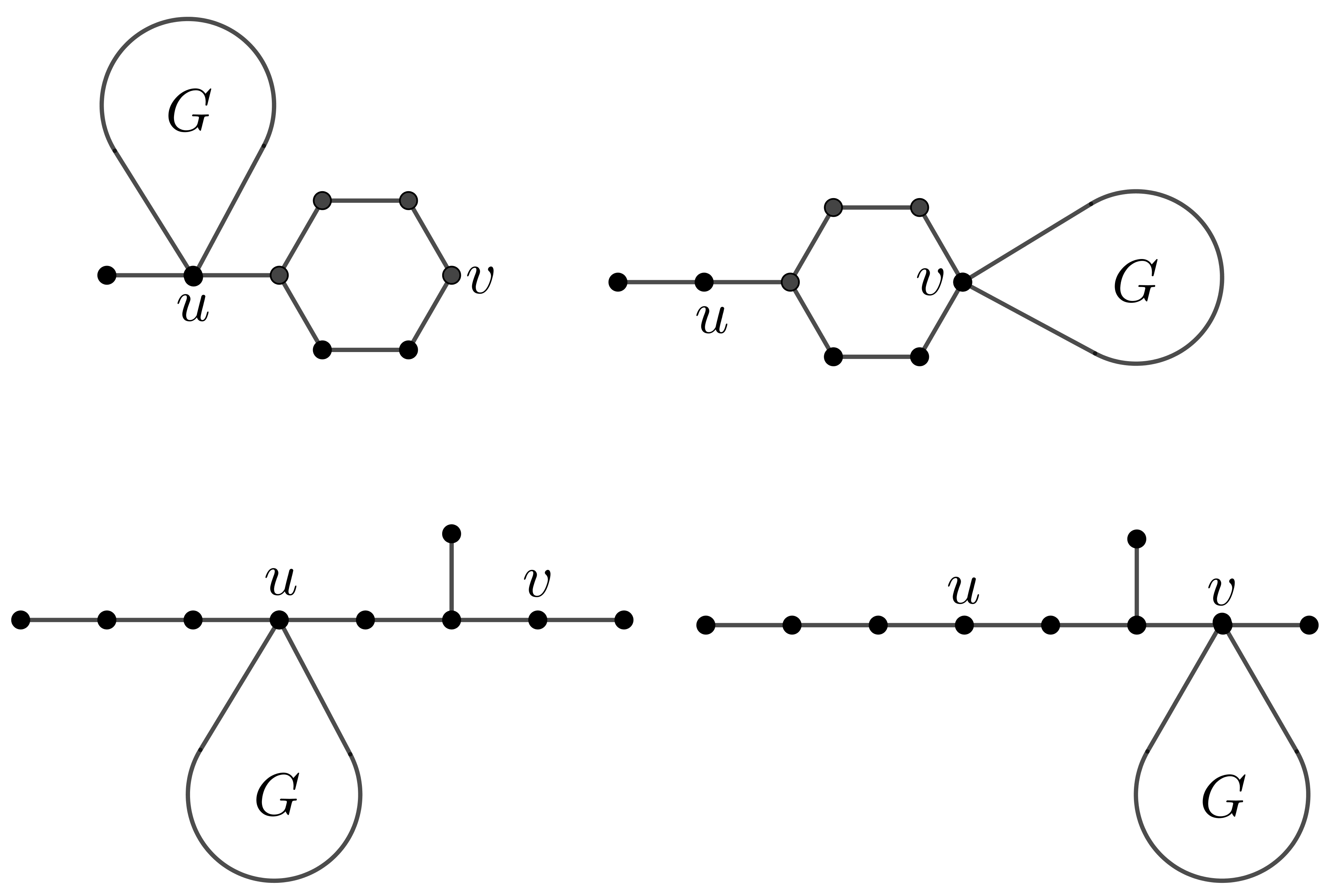}
\end{center}
\caption{Examples of pairs of cospectral graphs obtained 
         by attaching a rooted graph at isospectral vertices.}
\label{fig-2}
\end{figure}

Besides finding specific examples of pairs of cospectral graphs,
researchers in the 1970s were keen to provide unified approaches 
that allowed the construction of arbitrary numbers of new pairs of cospectral graphs.
This has led to further development of theoretical methods to investigate 
the characteristic polynomials of graphs and the corresponding walk generating functions.
Prominent examples related to the construction of arbitrary numbers of new pairs of cospectral graphs 
were given by Herndon and Ellzey~\cite{heel-1,heel-2},
Schwenk~\cite{schwenk-1,schwenk-2} and Godsil and McKay~\cite{gomc-2}.
Herndon and Ellzey~\cite{heel-1,heel-2} constructed several examples that are
based on the existence of pairs of isospectral vertices in a graph:
vertices $u$ and $v$ of the graph $G$ are isospectral
if the graphs $G-u$ and $G-v$ are cospectral, but not isomorphic.
Two examples of such vertex pairs in small graphs are shown in Fig.~\ref{fig-1}.
Using isospectral vertices we can obtain a pair of cospectral graphs simply 
by identifying either of the isospectral vertices with the root of an arbitrary rooted graph~$G$,
as indicated in Fig.~\ref{fig-2}.
The pair of graphs shown in the lower half of Fig.~\ref{fig-2}
was used by Schwenk to prove his celebrated result~\cite{schwenk-1}
that almost every tree is cospectral to another tree,
by showing that the proportion of trees 
that have either of the forms shown on the bottom of Fig.~\ref{fig-2}
tends to one as the number of their vertices tends to infinity.

Schwenk~\cite{schwenk-2} further generalized this coalescence approach as follows. 
For a graph~$G$ with a subset of vertices $S=\{s_1,\dots,s_k\}$,
a graph~$H$ with a subset of vertices $T=\{t_1,\dots,t_k\}$,
and a correspondence map $\theta\colon S\to T$ defined by $\theta(s_i)=t_i$,
we say that $S$ and $T$ are \emph{removal-cospectral sets} 
if for each $A\subseteq S$
the graphs $G-A$ and $H-\theta(A)$ are cospectral.
In such a case, 
we can take an arbitrary third graph~$J$ and a subset of its vertices $\{u_1,\dots,u_k\}$
to obtain a pair of cospectral graphs by forming two multiple coalescences:
one is obtained from $G$ and~$J$ by identifying vertices $s_i$ and~$u_i$ for $i=1,\dots,k$,
while the other is obtained from $H$ and~$J$ by identifying vertices $t_i$ and~$u_i$ for $i=1,\dots,k$.

Godsil and McKay~\cite{gomc-2} presented another celebrated approach to construct cospectral graphs,
based on which they were able to produce an estimated 72\% of the 51039 graphs on nine vertices 
that do not have unique spectrum.
For a graph~$G$, let $\pi=(C_1,\dots,C_k,D)$ be a partition of its vertex set
and suppose that for each $1\leq i,j\leq k$ and $v\in D$ the following conditions are satisfied:
a) any two vertices in $C_i$ have the same number of neighbours in $C_j$ and
b) $v$ has either $0$, $|C_i|/2$ or $|C_i|$ neighbours in~$C_i$.
The graph $G^{(\pi)}$ formed by local switching in~$G$ with respect to~$\pi$ is obtained from~$G$ as follows.
For each $v\in D$ and $1\leq i\leq k$ 
so that $v$ has $|C_i|/2$ neighbours in $C_i$,
the edges between $v$ and these $|C_i|/2$ neighbours are deleted,
while the edges between $v$ and the remaining $|C_i|/2$ vertices in~$C_i$ are added.
Then the graphs $G$ and $G^{(\pi)}$ are cospectral,
and moreover, their complements are also cospectral.
As a consequence, Godsil and McKay~\cite{gomc-2} improved Schwenk's result~\cite{schwenk-2}
by showing that almost all trees are cospectral with cospectral complements.
It should be noted, however, that 
Schwenk's removal-cospectral sets and Godsil-McKay local switching
are actually equivalent ways for constructing cospectral graphs~\cite{gomc-2}.

Several computational enumerations of cospectral graphs were presented in the literature.
The first exhaustive enumeration of such graphs for up to 9~vertices and trees for up to 14~vertices
was performed as early as 1976~\cite{gomc-1}.
The enumeration of cospectral graphs with 10 vertices was then presented by Lepovi\'c~\cite{lepo}.
Later, in a more recent work, that essentially revived interest in these topics, 
Haemers and  Spence~\cite{hasp} enumerated cospectral graphs with  11 vertices.
Namely, Haemers and Spence~\cite{hasp} observed that the percentage of cospectral graphs with 10 vertices is~21.3\%, 
which is slightly greater than that for 11 vertices (21.1\%).
Motivated by this observation 
(and the fact that almost all random matrices have simple eigenvalues~\cite{tavu},
which prevents the possibility of local switching to some extent)
van Dam and Haemers~\cite{daha-2} proposed the following conjecture:
almost all graphs ought to be determined by their spectrum.
As a result, considerable research has been conducted to establish spectral characterizations of various special types of graphs
(refer to~\cite{daha-1,daha-2} for relatively early findings).
The enumeration of cospectral graphs with 12 vertices presented by Brouwer and Spence~\cite{brsp}
showed that the percentage of cospectral graphs with 12 vertices further decreases to 18.8\%, 
providing additional support to van Dam and Haemers' expectations~\cite{daha-2}.

\begin{figure}[ht]
\begin{center}
\includegraphics[scale=0.75]{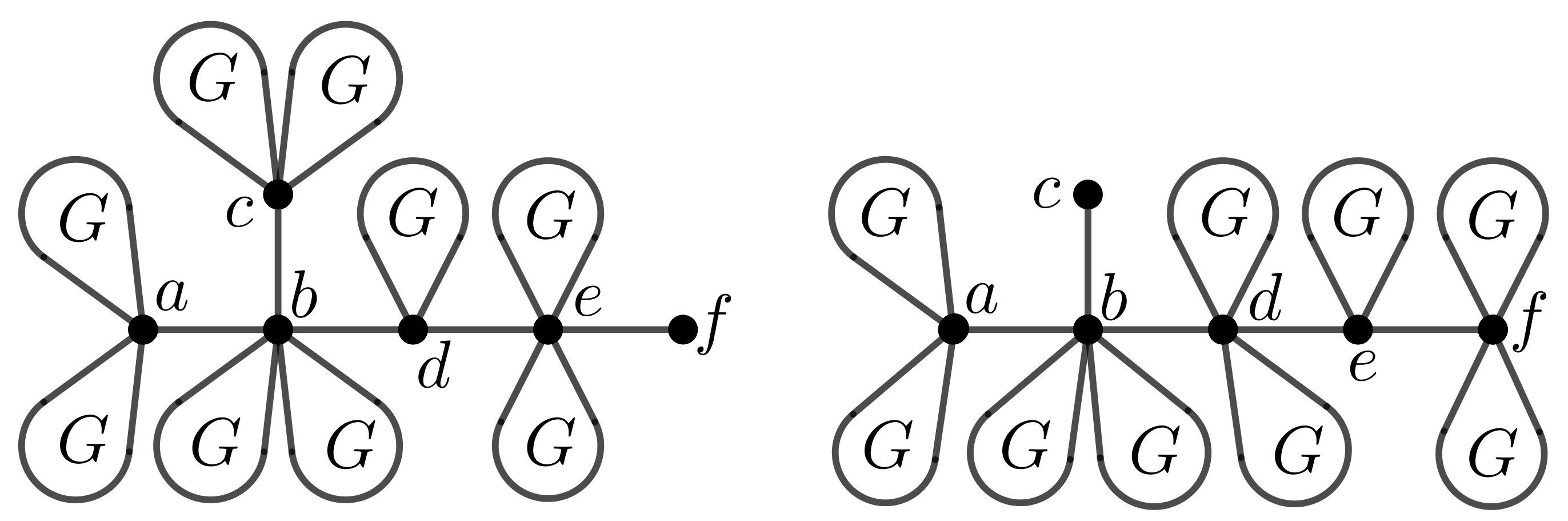}
\end{center}
\caption{A particular example of cospectral multiple coalescences from~\cite{hosoya-2}.}
\label{fig-3}
\end{figure}

More recently, the well-known theoretical chemist Haruo Hosoya 
drew the attention in a series of papers~\cite{hosoya-1,hosoya-2,hosoya-3}
to a particular aspect of constructing cospectral graphs by using coalescence of graphs:
besides attaching one copy of a rooted graph to either of two isospectral vertices,
pairs of cospectral graphs can be obtained by attaching multiple copies of a rooted graph
to different multisets of vertices in the underlying graph.
Hosoya managed to find a number of pairs such as the one shown in Fig.~\ref{fig-3},
where cospectral graphs are obtained by attaching up to nine copies of the rooted graph~$G$
to the vertices of the underlying graph~$T$ in two different ways.
Since Hosoya's construction scheme is a special case of Schwenk's multiple coalescences,
one will obviously obtain cospectral graphs whenever the appropriate vertex subsets in~$T$ form removal-cospectral sets.
However, Hosoya's attachment of the copies of the same rooted graph facilitates obtaining cospectral coalescences
with less stringent conditions for the corresponding vertex subsets in~$T$.
We describe these conditions 
through an explicit formula (\ref{eq-main}) for the characteristic polynomial of Hosoya's coalescences
that we obtain in Section~\ref{sc-inclusion-exclusion}.

In his works~\cite{hosoya-1,hosoya-2,hosoya-3}, 
Hosoya used the so-called $Z$-counting polynomial,
which is equivalent to the characteristic polynomial in the case of trees,
but not in the case of graphs that contain cycles.
Hosoya studied cospectrality of coalescences in which 
two or three copies of a rooted graph are attached to the vertices of an underlying tree,
and proposed a general expectation that
the characteristic polynomial of the multiple coalescence depends
on the family of characteristic polynomials of vertex-deleted subgraphs of the underlying tree~$T$.
We substantiate this expectation in Section~\ref{sc-inclusion-exclusion} 
through the explicit expression~(\ref{eq-main}) for the characteristic polynomial of Hosoya's coalescences.

Hosoya focused on the cases when the underlying tree is a path or contains a perfect matching and 
recommended identification of cospectral multiple coalescences in the cases when the underlying graph contains cycles.
In Section~\ref{sc-enumeration} we computationally enumerate such multiple coalescences
in the cases when the underlying graph is
a unicyclic graph, a small catacondensed benzenoid or a general graph,
in addition to furthering enumeration in cases when the underlying graph is a path, a tree or a tree with a perfect matching.
At the end, we identify an infinite family of cospectral multiple coalescences with paths as underlying graphs,
whose existence was deemed rare by Hosoya.

\section{Definitions and preliminaries}

We first introduce necessary notation for multiple coalescences with the same underlying graph.

\begin{definition}
Let $T$ and $G_1,\dots,G_k$ be vertex-disjoint graphs for some $k\geq 1$.
Let $u_1,\dots,u_k$ be the distinct vertices of~$T$,
and for each $i=1,\dots,k$, 
let $v_i$ be a vertex of the graph $G_i$.
The (multiple) coalescence denoted by 
$$
T(u_1\=v_1)G_1\cdots(u_k\=v_k)G_k
$$
is the graph obtained from the union $T\cup G_1\cup\cdots\cup G_k$
by identifying vertices $u_i$ and $v_i$ for each $i=1,\dots,k$.
\end{definition}

In his papers~\cite{hosoya-1,hosoya-2,hosoya-3},
Hosoya was mostly interested in the case 
when $G_1,\dots,G_k$ are copies of the same rooted graph~$G$ with the root~$r$,
and in such a way that several copies of $G$ may be attached at a single vertex of~$T$.
To accommodate this setup, let us introduce further notation.

\begin{definition}
For a rooted graph $G$ with the root~$r$ and $a\geq 1$,
let $G^{(a)}$ denote the coalescence $G(r\=r)G\cdots(r\=r)G$
in which the roots of $a$ copies of~$G$ are all mutually identified,
and also denoted by~$r$.
\end{definition}

From the above definition,
we trivially have that $G^{(1)}\cong G$ and 
that $G^{(a)}\-r$ is the disjoint union of $a$ copies of~$G\-r$.

Hosoya's examples of multiple coalescences may now be described in terms of the underlying graph~$T$,
distinct vertices $(u_1,\dots,u_k)$ of $T$ at which copies of $G$ are attached, and
the {\em signature} $(a_1,\dots,a_k)$ stating that $a_i$ copies of~$G$ are attached at~$u_i$ for $i=1,\dots,k$.
For example, two coalescences in Fig.~\ref{fig-3} may be denoted as 
$$
T(a\=r)G^{(2)}(b\=r)G^{(2)}(c\=r)G^{(2)}(e\=r)G^{(2)}(d\=r)G^{(1)}
$$
and
$$
T(a\=r)G^{(2)}(b\=r)G^{(2)}(d\=r)G^{(2)}(f\=r)G^{(2)}(e\=r)G^{(1)},
$$
respectively, where $T$ is the underlying tree with the vertex set $\{a,b,c,d,e,f\}$
and the signature in both coalescences is $(2,2,2,2,1)$.

Let us denote by $P(G)$, or simply by $PG$ when there is no confusion,
the characteristic polynomial of adjacency matrix of the graph~$G$ in terms of the variable~$x$.
Schwenk~\cite{schwenk-1} proved the following formula for the characteristic polynomial
of the coalescence of two graphs:
\begin{equation}
\label{eq-schwenk}
PG(u\=v)H = PG P(H\-v) + P(G\-u)PH - x P(G\-u)P(H\-v),
\end{equation}
where $G\-u$ denotes deletion of the vertex~$u$ and its incident edges from the graph~$G$
(and likewise for $H\-v$).

A simple inductive argument yields the characteristic polynomial of $G^{(a)}$.
\begin{lemma}
\label{le-Gpower}
If $G$ is a rooted graph with the root~$r$, then for any $a\geq 1$
\begin{equation}
\label{eq-Gpower}
PG^{(a)}=[aPG  - (a\-1)xP(G\-r)]P(G\-r)^{a-1}.
\end{equation}
\end{lemma}
\begin{proof}
This is trivially satisfied for $a=1$. Assume thus that the statement holds for some $a\geq 1$.
From $G^{(a+1)}\cong G(r=r)G^{(a)}$ and 
the fact that $G^{(a)}-r$ is the union of $a$~disjoint copies of $G\!-\!r$,
we have
\begin{align*}
PG^{(a+1)} 
  &= PG P(G^{(a)}\-r) + P(G\-r)PG^{(a)} - xP(G\-r)P(G^{(a)}\-r) \\
  &= PG P(G\-r)^a + P(G\-r)[a PG - (a\-1)xP(G\-r)]P(G\-r)^{a-1} - xP(G\-r)P(G\-r)^a \\
  &= [(a\+1) PG - axP(G\-r)]P(G\-r)^a.
\end{align*}
\end{proof}

\section{Characteristic polynomial of Hosoya's coalescences}
\label{sc-inclusion-exclusion}

Here we state our main theorem that brings 
an explicit inclusion-exclusion style formula for the characteristic polynomial of multiple coalescences of Hosoya type.

\begin{theorem}
Let $T$ be a graph with selected distinct vertices $u_1,\dots,u_k$ for some $k\geq 1$.
For each $i=1,\dots,k$ let $G_i$ be a rooted graph with the root~$r_i$
and let 
$$
Q_i = PG_i - xP(G_i\-r_i)
\quad\mbox{and}\quad
R_i = P(G_i\-r_i).
$$
Then for any signature $(a_1,\dots, a_k)$, we have 
\begin{equation}
\label{eq-main}
PT(u_1\=r_1)G_1^{(a_1)}\cdots(u_k\=r_k)G_k^{(a_k)} 
 = \sum_{I\subseteq\{1,\dots,k\}} 
    P(T\-\sum_{i\in I}u_i) \prod_{j\in I} a_j \prod_{l\in I} Q_l \prod_{m=1}^k R_m^{a_m-\left|\{m\}\cap I\right|}.
\end{equation}
\end{theorem}

\begin{proof}
Note that for $k=0$ the formula~(\ref{eq-main}) reduces to the obvious identity $PT=PT$.
However, we need the case $k=1$ in the proof of the inductive step,
so we take $k=1$ as the basis of induction.

For $k=1$ the value of $I$ in the first sum in~(\ref{eq-main}) is either $\emptyset$ or $\{1\}$,
so that~(\ref{eq-main}) reduces to
$$
PT(u_1\=r_1)G_1^{(a_1)} = PT R_1^{a_1} + P(T\-u_1) a_1 Q_1 R_1^{a_1-1}.
$$
After replacing $Q_1$ and~$R_1$, from Lemma~\ref{le-Gpower} we have
\begin{align*}
PT(u_1\=r_1)G_1^{(a_1)}
 &= PT P(G_1\-r_1)^{a_1} + P(T\-u_1) a_1 [PG_1 - xP(G_1\-r_1)] P(G_1\-r_1)^{a_1-1} \\
 &= PT P(G_1^{(a_1)}\-r_1) + P(T\-u_1) PG_1^{(a_1)} - x P(T\-u_1)P(G_1^{(a_1)}\-r_1),
\end{align*}
which is correct by the Schwenk's formula~(\ref{eq-schwenk}).

Assume now that~(\ref{eq-main}) holds for some $k\geq 1$. 
The multiple coalescence 
$$
PT(u_1\=r_1)G_1^{(a_1)}\cdots(u_k\=r_k)G_k^{(a_k)}(u_{k+1}\=r_{k+1})G_{k+1}^{(a_{k+1})}
$$ 
is, at the same time, also a coalescence of $PT(u_1\=r_1)G_1^{(a_1)}\cdots(u_k\=r_k)G_k^{(a_k)}$ with $G_{k+1}^{(a_{k+1})}$.
From the basis of induction we then have 
\begin{align*}
  & PT(u_1\=r_1)G_1^{(a_1)}\cdots(u_k\=r_k)G_k^{(a_k)}(u_{k+1}\=r_{k+1})G_{k+1}^{(a_{k+1})} \\
  &\quad = PT(u_1\=r_1)G_1^{(a_1)}\cdots(u_k\=r_k)G_k^{(a_k)} R_{k+1}^{a_{k+1}}
       \qquad\hspace{97pt}\mbox{(for $I=\emptyset$)} \\
  &\quad + P(T\-u_{k+1})(u_1\=r_1)G_1^{(a_1)}\cdots(u_k\=r_k)G_k^{(a_k)} a_{k+1}Q_{k+1}R_{k+1}^{a_{k+1}-1}
       \qquad\mbox{(for $I=\{k+1\}$)},
\end{align*}
while from the inductive assumption we further obtain
\begin{align*}
  & PT(u_1\=r_1)G_1^{(a_1)}\cdots(u_k\=r_k)G_k^{(a_k)}(u_{k+1}\=r_{k+1})G_{k+1}^{(a_{k+1})} \\
  &\quad = \sum_{I\subseteq\{1,\dots,k\}} 
                 P(T\-\sum_{i\in I}u_i) \prod_{j\in I} a_j \prod_{l\in I} Q_l \prod_{m=1}^k R_m^{a_m-\left|\{m\}\cap I\right|}
                 \cdot R_{k+1}^{a_{k+1}} \\
  &\quad + \sum_{I\subseteq\{1,\dots,k\}} 
                 P(T\-u_{k+1}\-\sum_{i\in I}u_i) \prod_{j\in I} a_j\cdot a_{k+1} \prod_{l\in I} Q_l\cdot Q_{k+1} 
                 \prod_{m=1}^k R_m^{a_m-\left|\{m\}\cap I\right|}\cdot R_{k+1}^{a_{k+1}-1}.
\end{align*}
Setting $I^\ast\=I$ in the first and $I^\ast\=I\cup\{k+1\}$ in the second sum over~$I$ above we get
\begin{align*}
  & PT(u_1\=r_1)G_1^{(a_1)}\cdots(u_k\=r_k)G_k^{(a_k)}(u_{k+1}\=r_{k+1})G_{k+1}^{(a_{k+1})} \\
  &\quad = \sum_{\substack{I^\ast\subseteq\{1,\dots,k,k+1\}\\k+1\notin I^\ast}} 
                 P(T\-\sum_{i\in I^\ast}u_i) \prod_{j\in I^\ast} a_j \prod_{l\in I^\ast} Q_l 
                                                          \prod_{m=1}^{k+1} R_m^{a_m-\left|\{m\}\cap I^\ast\right|} \\
  &\quad + \sum_{\substack{I^\ast\subseteq\{1,\dots,k,k+1\}\\k+1\in I^\ast}} 
                 P(T\-\sum_{i\in I^\ast}u_i) \prod_{j\in I^\ast} a_j \prod_{l\in I^\ast} Q_l 
                                                          \prod_{m=1}^{k+1} R_m^{a_m-\left|\{m\}\cap I^\ast\right|} \\
  &\quad = \sum_{I^\ast\subseteq\{1,\dots,k,k+1\}} 
                 P(T\-\sum_{i\in I^\ast}u_i) \prod_{j\in I^\ast} a_j \prod_{l\in I^\ast} Q_l 
                                                          \prod_{m=1}^{k+1} R_m^{a_m-\left|\{m\}\cap I^\ast\right|},
\end{align*}
which proves the inductive step.
\end{proof}

In the case when all rooted graphs $G_1,\dots,G_k$ are equal to~$G$
and all their roots are equal to~$r$, we obtain the following corollary.
\begin{corollary}
Let $T$ be a graph with selected distinct vertices $u_1,\dots,u_k$ for some $k\geq 1$.
For a rooted graph $G$ with the root~$r$
let 
$$
Q = PG - xP(G\-r)
\quad\mbox{and}\quad
R = P(G\-r).
$$
Then for any signature $(a_1,\dots, a_k)$, we have 
\begin{equation}
\label{eq-main-2}
PT(u_1\=r)G^{(a_1)}\cdots(u_k\=r)G^{(a_k)} 
 = \sum_{I\subseteq\{1,\dots,k\}} 
    P(T\-\sum_{i\in I}u_i) \prod_{j\in I} a_j\, Q^{|I|} R^{\sum_{m=1}^k a_m-|I|} .
\end{equation}
\end{corollary}

\begin{example}
\label{ex-1}
Let us use (\ref{eq-main-2}) to obtain the characteristic polynomials of two multiple coalescences shown in Fig.~\ref{fig-3}.
Let $T$ be the underlying tree consisting of vertices $\{a,b,c,d,e,f\}$.
Grouping the subsets $I\subseteq\{1,\dots,5\}$ by their cardinality,
we obtain

\vspace{-6pt}
{\small
\begin{align*}
& PT(a\=r)G^{(2)}(b\=r)G^{(2)}(c\=r)G^{(2)}(e\=r)G^{(2)}(d\=r)G^{(1)} \\
& = PT R^9 \\
& + [2P(T\-a) \+ 2P(T\-b) \+ 2P(T\-c) \+ 2P(T\-e) \+ P(T\-d)]QR^8 \\
& + [4P(T\-a\-b)\+4P(T\-a\-c)\+4P(T\-a\-e)\+2P(T\-a\-d)\+4P(T\-b\-c) \\
& \phantom{+} \+4P(T\-b\-e)\+2P(T\-b\-d)\+4P(T\-c\-e)\+2P(T\-c\-d)\+2P(T\-e\-d)]Q^2R^7 \\
& + [8P(T\-a\-b\-c)\+8P(T\-a\-b\-e)\+4P(T\-a\-b\-d)\+8P(T\-a\-c\-e)\+4P(T\-a\-c\-d) \\
& \phantom{+} \+4P(T\-a\-e\-d)\+8P(T\-b\-c\-e)\+4P(T\-b\-c\-d)\+4P(T\-b\-e\-d)\+4P(T\-c\-e\-d)]Q^3R^6 \\
& + [16P(T\-a\-b\-c\-e)\+8P(T\-a\-b\-c\-d)\+8P(T\-a\-b\-e\-d) \\
& \phantom{+} \+8P(T\-a\-c\-e\-d)\+8P(T\-b\-c\-e\-d)]Q^4R^5 \\
& + 16P(T\-a\-b\-c\-e\-d)Q^5R^4 \\
& = PT R^9 \+ (9x^5\-29x^3 + 14x)QR^8 \+ (32x^4\-58x^2\+8)Q^2R^7 \+ (56x^3 \- 44x)Q^3R^6 \+ (48x^2\-8)Q^4R^5 \+ 16x Q^5R^4.
\end{align*}
}

\vspace{-6pt}
\noindent
and

\vspace{-6pt}
{\small
\begin{align*}
& PT(a\=r)G^{(2)}(b\=r)G^{(2)}(d\=r)G^{(2)}(f\=r)G^{(2)}(e\=r)G^{(1)} \\
& = PT R^9 \\
& + [2P(T\-a) \+ 2P(T\-b) \+ 2P(T\-d) \+ 2P(T\-f) \+ P(T\-e)]QR^8 \\
& + [4P(T\-a\-b)\+4P(T\-a\-d)\+4P(T\-a\-f)\+2P(T\-a\-e)\+4P(T\-b\-d) \\
& \phantom{+} \+4P(T\-b\-f)\+2P(T\-b\-e)\+4P(T\-d\-f)\+2P(T\-d\-e)\+2P(T\-f\-e)]Q^2R^7 \\
& + [8P(T\-a\-b\-d)\+8P(T\-a\-b\-f)\+4P(T\-a\-b\-e)\+8P(T\-a\-d\-f)\+4P(T\-a\-d\-e) \\
& \phantom{+} \+4P(T\-a\-f\-e)\+8P(T\-b\-d\-f)\+4P(T\-b\-d\-e)\+4P(T\-b\-f\-e)\+4P(T\-d\-f\-e)]Q^3R^6 \\
& + [16P(T\-a\-b\-d\-f)\+8P(T\-a\-b\-d\-e)\+8P(T\-a\-b\-f\-e) \\
&  \phantom{+} \+8P(T\-a\-d\-f\-e)\+8P(T\-b\-d\-f\-e)]Q^4R^5 \\
& + 16P(T\-a\-b\-d\-f\-e)Q^5R^4 \\
& = PT R^9 \+ (9x^5\-29x^3 \+ 14x)QR^8 \+ (32x^4\-58x^2\+8)Q^2R^7 \+ (56x^3\-44x)Q^3R^6 \+ (48x^2\-8)Q^4R^5 \+ 16xQ^5R^4,
\end{align*}
}

\vspace{-6pt}
\noindent
showing that the two coalescences are cospectral regardless of the rooted graph~$G$.

Note, however, that the vertex subsets $\{a,b,c,e,d\}$ and $\{a,b,d,f,e\}$,
at which the copies of $G$ are attached, are not removal-cospectral.
Namely, apart from $a$ and $c$, all vertex-deleted subgraphs of~$T$ have distinct characteristic polynomials:
\begin{align*}
P(T\-a)\=P(T\-c)&\=x^5\-4x^3\+3x, \\
P(T\-b)&\=x^5-2x^3, \\
P(T\-d)&\=x^5-3x^3+2x, \\
P(T\-e)&\=x^5-3x^3, \\
P(T\-f)&\=x^5-4x^3+2x.
\end{align*}
Hence, there is no way to establish a one-to-one correspondence $\theta\colon\{a,b,c,e,d\}\to\{a,b,d,f,e\}$
such that both $T-\theta(a)$ and $T-\theta(c)$ are cospectral (to the cospectral subgraphs $T-a$ and~$T-c$).
\end{example}

We can now discuss the consequences of expressions (\ref{eq-main}) and~(\ref{eq-main-2}).
First, these expressions substantiate Hosoya's expectations 
about the $Z$-counting polynomial from~\cite[Conjecture~1]{hosoya-3},
as the $Z$-counting polynomial is identical to the characteristic polynomial 
when both the underlying graph $T$ and the attached rooted graphs $G_1,\dots,G_k$ are trees.
One should note here, however, that Conjecture~1 in~\cite{hosoya-3} is stated rather informally,
so that the expressions (\ref{eq-main}) and~(\ref{eq-main-2}) 
should not be considered directly as the proof of Hosoya's conjecture,
but rather as its clarification and formalisation.

Next we clarify conditions under which certain multiple coalescences are cospectral.
\begin{corollary}
\label{co-removal-cospectral}
Let $T_1$ be a graph with selected distinct vertices $u_1,\dots,u_k$,
and $T_2$ a graph with selected distinct vertices $v_1,\dots,v_k$ for some $k\geq 1$.
For a fixed signature $(a_1,\dots,a_k)$,
the multiple coalescences
\begin{equation}
\label{eq-free-coalescences}
T_1(u_1\=r_1)G_1^{(a_1)}\cdots(u_k\=r_k)G_k^{(a_k)}
\quad\mbox{and}\quad
T_2(v_1\=r_1)G_1^{(a_1)}\cdots(v_k\=r_k)G_k^{(a_k)}
\end{equation}
are cospectral for all possible choices of the rooted graphs $G_1,\dots,G_k$ and their roots $r_1,\dots,r_k$
if and only if
$T_1$ and $T_2$ are cospectral graphs 
with the removal-cospectral sets $\{u_1,\dots,u_k\}$ and $\{v_1,\dots,v_k\}$.
\end{corollary}

\begin{proof}
In one direction,
if $T_1$ and $T_2$ are cospectral graphs 
with removal-cospectral sets $\{u_1,\dots,u_k\}$ and $\{v_1,\dots,v_k\}$,
then by definition
$$
P(T_1\-\sum_{i\in I} u_i)\=P(T_2-\sum_{i\in I} v_i)
$$
for each $I\subseteq\{1,\dots,k\}$ and
the cospectrality of the coalescences (\ref{eq-free-coalescences}) follows directly from Eq.~(\ref{eq-main}).

On the other hand, 
one can always proclaim $G_i^{(a_i)}$ to be a new rooted graph $G^{\ast}_i$,
so that cospectrality of the coalescences~(\ref{eq-free-coalescences}) must also hold for the signature $(1,\dots,1)$.
As the rooted graphs $G_1,\dots,G_k$ with roots $r_1,\dots,r_k$ are chosen independently,
one has to assume that the polynomials $Q_1,R_1,\dots,Q_k,R_k$ are mutually independent.
From this independence and Eq.~(\ref{eq-main}) for the signature $(1,\dots,1)$,
we immediately conclude that the coalescences in~(\ref{eq-free-coalescences}) are cospectral 
for arbitrary rooted graphs $G_1,\dots,G_k$ with roots $r_1,\dots,r_k$
if and only if
$$
P(T_1\-\sum_{i\in I} u_i)\=P(T_2-\sum_{i\in I} v_i)
$$
for each $I\subseteq\{1,\dots,k\}$,
i.e., if and only if $T_1$ and $T_2$ are cospectral graphs 
with the removal-cospectral sets $\{u_1,\dots,u_k\}$ and $\{v_1,\dots,v_k\}$.
\end{proof}

The situation is, however, different in the case that Hosoya considered in~\cite{hosoya-1,hosoya-2,hosoya-3}.
\begin{corollary}
\label{co-Hosoya}
Let $T_1$ be a graph with selected distinct vertices $u_1,\dots,u_k$,
and $T_2$ a graph with selected distinct vertices $v_1,\dots,v_k$ for some $k\geq 1$.
For fixed signatures $(a_1,\dots,a_k)$ and $(b_1,\dots,b_k)$
satisfying $a_1\geq\cdots\geq a_k$ and $b_1\geq\cdots\geq b_k$,
the multiple coalescences
\begin{equation}
\label{eq-bound-coalescences}
T_1(u_1\=r)G^{(a_1)}\cdots(u_k\=r)G^{(a_k)}
\quad\mbox{and}\quad
T_2(v_1\=r)G^{(b_1)}\cdots(v_k\=r)G^{(b_k)}
\end{equation}
are cospectral for all possible choices of the rooted graph~$G$ and its root~$r$
if and only if $(a_1,\dots,a_k)=(b_1,\dots,b_k)$ and
\begin{equation}
\label{eq-Hosoya-condition}
\sum_{\substack{I\subseteq\{1,\dots,k\}\\|I|=l}} P(T_1\-\sum_{i\in I}u_i)\prod_{j\in I}a_j
=\sum_{\substack{I\subseteq\{1,\dots,k\}\\|I|=l}} P(T_2\-\sum_{i\in I}v_i)\prod_{j\in I}a_j
\end{equation}
for each $0\leq l\leq k$.
\end{corollary}

\begin{proof}
In the coalescences~(\ref{eq-bound-coalescences})
one can consider the polynomials $Q=PG-xP(G-r)$ and $R=P(G-r)$ to be mutually independent,
as one can arbitrarily choose the subgraph~$G-r$ and then arbitrarily connect the root~$r$ with the vertices of~$G-r$.
In such case 
Eq.~(\ref{eq-main-2}) implies that the coalescences~(\ref{eq-bound-coalescences}) are cospectral
if and only if 
$
\sum_{m=1}^k a_m = \sum_{m=1}^k b_m
$
and
\begin{equation}
\label{eq-Hosoya-condition-again}
\sum_{\substack{I\subseteq\{1,\dots,k\}\\|I|=l}} P(T_1\-\sum_{i\in I}u_i)\prod_{j\in I}a_j
=\sum_{\substack{I\subseteq\{1,\dots,k\}\\|I|=l}} P(T_2\-\sum_{i\in I}v_i)\prod_{j\in I}b_j
\end{equation}
for each $0\leq l\leq k$.
For each $I\subseteq\{1,\dots,k\}$ with $|I|=l$,
the characteristic polynomials $P(T_1-\sum_{i\in I} u_i)$ and $P(T_2-\sum_{i\in I}v_i)$ are 
monic of degrees $n_1-l$ and $n_2-l$, respectively,
where $n_1$ and~$n_2$ denote the numbers of vertices in $T_1$ and~$T_2$, respectively.
Hence for each $0\leq l\leq k$,
the left-hand side sum in~(\ref{eq-Hosoya-condition-again}) is
a polynomial of degree $n_1-l$ with the leading coefficient 
$\sum_{I\subseteq\{1,\dots,k\}, |I|=l} \prod_{j\in I}a_j$,
while the right-hand side sum is
a polynomial of degree $n_2-l$ with the leading coefficient
$\sum_{I\subseteq\{1,\dots,k\}, |I|=l} \prod_{j\in I}b_j$.
Thus $n_1=n_2$ and
$$
\sum_{I\subseteq\{1,\dots,k\}, |I|=l} \prod_{j\in I}a_j
=\sum_{I\subseteq\{1,\dots,k\}, |I|=l} \prod_{j\in I}b_j
$$
for each $0\leq l\leq k$.
As this means that 
each elementary symmetric polynomial has equal values 
for the variables $a_1,\dots,a_k$ and the variables $b_1,\dots,b_l$,
we conclude that the families $a_1,\dots,a_k$ and $b_1,\dots,b_k$ are equal,
i.e., that $(a_1,\dots,a_k)=(b_1,\dots,b_k)$ due to $a_1\geq\cdots\geq a_k$ and $b_1\geq\cdots\geq b_k$.
From this, the conditions~(\ref{eq-Hosoya-condition-again}) directly translate to the conditions~(\ref{eq-Hosoya-condition}).

The other direction is trivial: 
if $(a_1,\dots,a_k)=(b_1,\dots,b_k)$ and 
the condition~(\ref{eq-Hosoya-condition}) holds for each $0\leq l\leq k$,
then the multiple coalescences~(\ref{eq-bound-coalescences}) are cospectral by Eq.~(\ref{eq-main-2}).
\end{proof}

Due to the more specialized structure of Hosoya's coalescences, 
we can see from Corollary~\ref{co-Hosoya} that
the subsets $\{u_1,\dots,u_k\}$ and $\{v_1,\dots,v_k\}$ need not be removal-cospectral 
in order to obtain cospectral coalescences.
This was the case in Example~\ref{ex-1},
in whose calculations a particular appearance of the condition~(\ref{eq-Hosoya-condition}) can be easily observed.

\section{Computational results}
\label{sc-enumeration}

Hosoya has found a number of examples of cospectral multiple coalescences (CMC) 
in his papers~\cite{hosoya-1,hosoya-2,hosoya-3} using back-of-envelope calculations. 
Here we have instead developed an extensive suite of Java classes to exhaustively search for examples of multiple coalescences,
based upon the existing Java framework for working with graphs~\cite{graph6java},
in which we have added an implementation of Samuelson-Berkowitz algorithm for computing characteristic polynomials~\cite{berk}.
These Java classes are available in source code at~\cite{zenodo},
together with two precompiled jar archives that may be run from a command line.
They may be used as follows:
\begin{itemize}
\item Suppose that we want to find examples of multiple coalescences
         in which the underlying graph belongs to a given set of graphs,
         which are collected in {\tt file.g6} in graph6 format~\cite{graph6java}, one graph per line.
\item We first compute the coefficients of characteristic polynomials of these graphs
         by issuing the command
         \begin{quote}
         {\tt java -jar listcharpolys.jar file.g6}
         \end{quote}
         in the terminal.
         This produces {\tt file.g6.charpoly.g6} which, together with a g6 code of each graph,
         contains a list of the coefficients of its characteristic polynomial in the same line.
         Note that {\tt listcharpolys.jar} is obtained by creating a jar archive
         from the {\tt main} method in the class {\tt ReporterTemplate.java} \cite{zenodo}.
\item Next we sort {\tt file.g6.charpoly.g6} to group cospectral graphs together.
         In Unix-based operating systems this may be done by issuing the terminal command
         \begin{quote}
         {\tt sort -n -k2 -o file.charpoly.sorted.g6 file.g6.charpoly.g6}
         \end{quote}
\item Finally, we start the search for examples of CMCs by issuing the command
         \begin{quote}
         {\tt java -jar hosoyacospectrality.jar file.charpoly.sorted.g6 <mse>}
         \end{quote}         
         in the terminal. Here {\tt <mse>} denotes the maximum possible signature entry (MSE).
         This command will process all groups of cospectral graphs from {\tt file.charpoly.sorted.g6},
         for each such group it will process all distinct signatures with entries between 1 and {\tt <mse>},
         and for each signature it will process all variations of the appropriate number of vertices from each cospectral graph,
         looking to identify examples that satisfy all conditions~(\ref{eq-Hosoya-condition}).
         The program will skip removal-cospectral sets of vertices, 
         which trivially yield CMCs by Corollary~\ref{co-removal-cospectral}.
         The program also requires that the signature entries are relatively prime,
         as for $d>1$
         $$
         T(u_1=r)G^{(da_1)}\cdots(u_k=r)G^{(da_k)}
         \cong T(u_1=r)\left(G^{(d)}\right)^{(a_1)}\cdots(u_k=r)\left(G^{(d)}\right)^{(a_k)},
         $$
         implying that the coalescences with the signature $(da_1,\dots,da_k)$ are 
         simply special cases of coalescences with the signature $(a_1,\dots,a_k)$.
         Note that {\tt hosoyacospectrality.jar} is obtained by creating a jar archive 
         from the {\tt main} method in the class {\tt HosoyaCospectrality.java} \cite{zenodo}.         
\end{itemize}

This exhaustive search quickly succumbs to combinatorial explosion, 
so that it can hardly be completed for graphs with more than a handful (10-12) of vertices,
regardless of the MSE.
During the execution of {\tt hosoyacospectrality.jar},
the examples of CMCs are saved to disk as soon as they are found,
so that even incomplete searches may still yield useful data.
The examples found are saved in Graphviz's dot format (see \url{graphviz.org}),
so that the information about the underlying graph, the signature entries and the selected vertices 
may be visualised by one of the Graphviz's layout programs, such as {\tt neato}.

\begin{figure}
\centering
\includegraphics[scale=0.65]{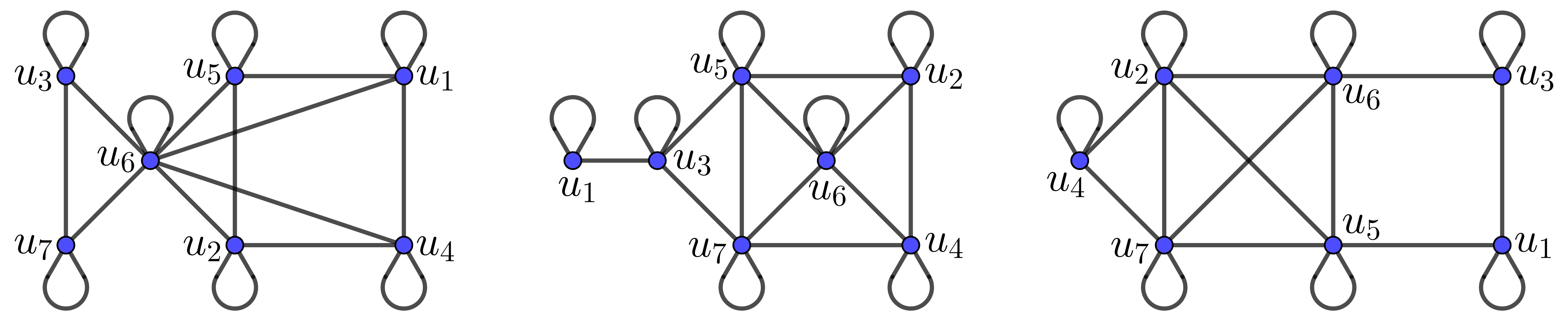}
{\small
\caption{A triplet of CMCs with underlying graphs on 7 vertices and the signature $(1,1,1,1,1,1,1)$.
In this and the subsequent figures, 
``droplets'' represent copies of the rooted graph that is attached at the corresponding vertices.
For each of these graphs and the vertices as labeled, 
the values of $\sum_{I\subseteq\{1,\dots,k\}, |I|=l} P(T_1-\sum_{i\in I}u_i)\prod_{j\in I}a_j$
for $l=0,\dots,7$ are as follows:
$y^{7} - 11y^{5} - 10y^{4} + 16y^{3} + 16y^{2}$, 
$7y^{6} - 55y^{4} - 40y^{3} + 48y^{2} + 32y$, 
$21y^{5} - 110y^{3} - 60y^{2} + 48y + 16$, 
$35y^{4} - 110y^{2} - 40y + 16$, 
$35y^{3} - 55y - 10$, 
$21y^{2} - 11$, 
$7y$, 
$1$.
\label{fig-example-1}
}}
\end{figure}

\begin{figure}
\centering
\includegraphics[scale=0.65]{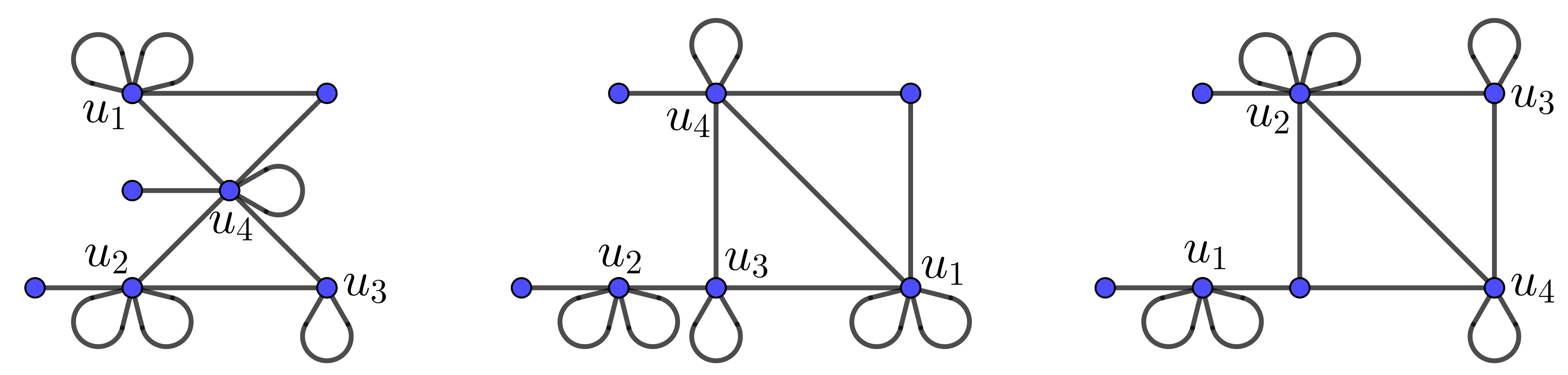} \\[6pt]
$a)$ \\[6pt]
\includegraphics[scale=0.65]{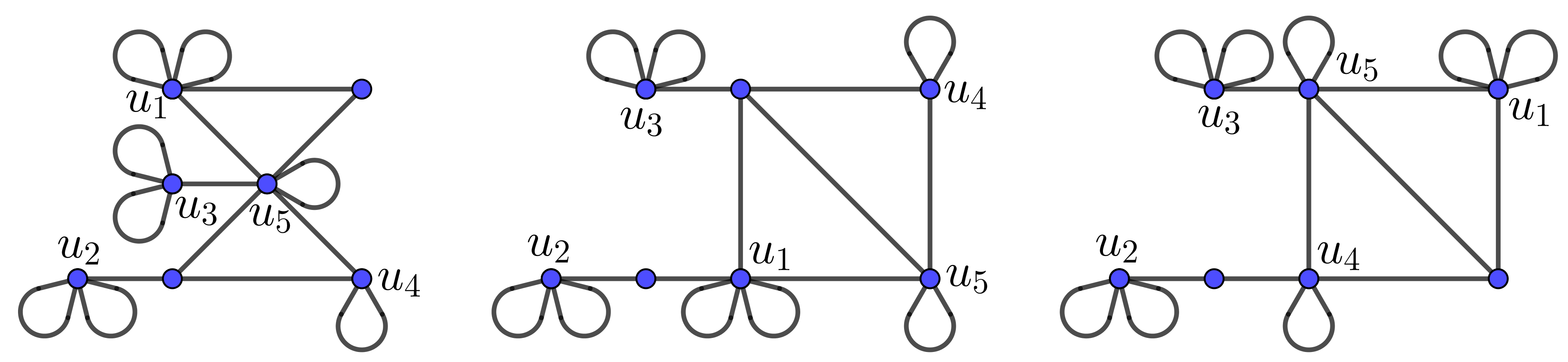} \\[6pt]
$b)$
{\small
\caption{Two triplets of CMCs with underlying graphs on 7 vertices and:
a) the signature $(2,2,1,1)$,
b) the signature $(2,2,2,1,1)$.
While the second and the third underlying graph in these coalescences are isomorphic,
the selections of vertices at which the copies of the rooted graph are attached
differ in these underlying graphs.
\label{fig-b}
}}
\end{figure}

\begin{figure}
\centering
\includegraphics[scale=0.65]{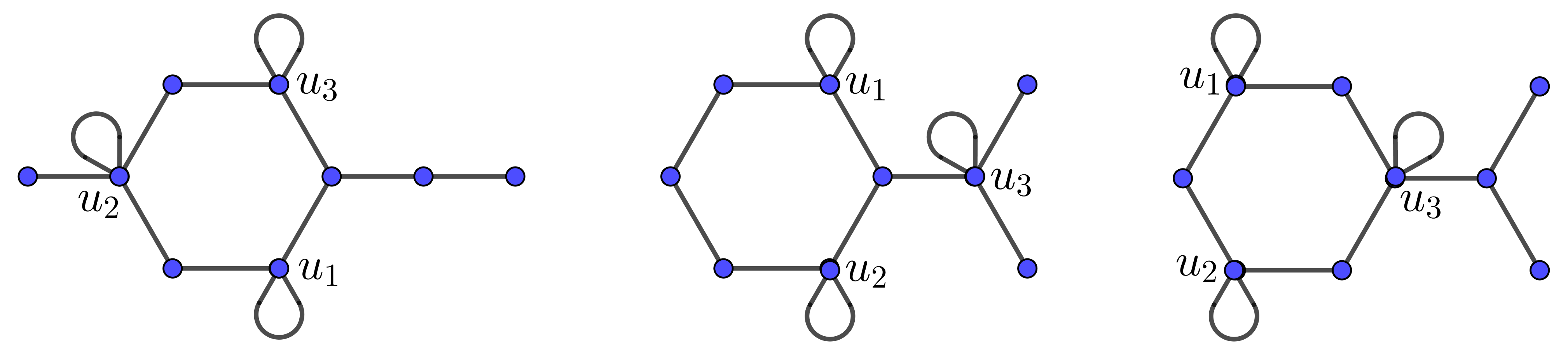} \\[6pt]
$a)$ \\[6pt]
\includegraphics[scale=0.65]{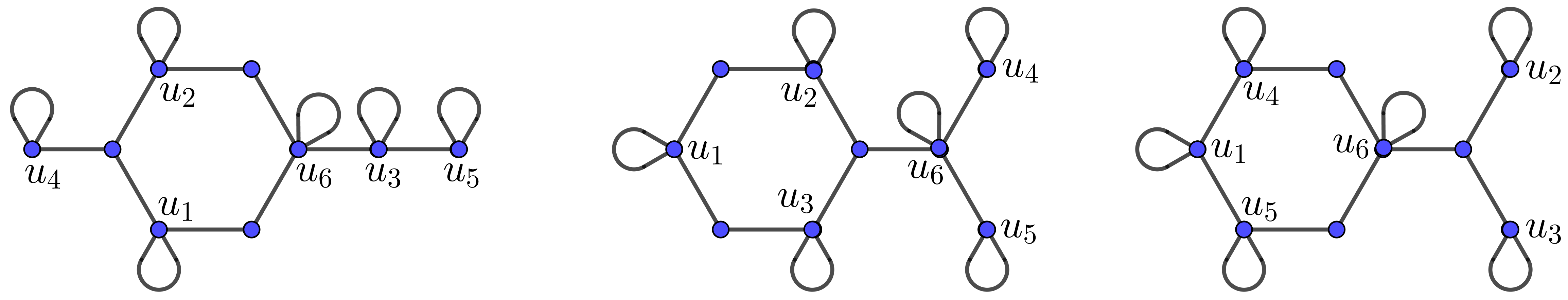} \\[6pt]
$b)$ \\[6pt]
\includegraphics[scale=0.65]{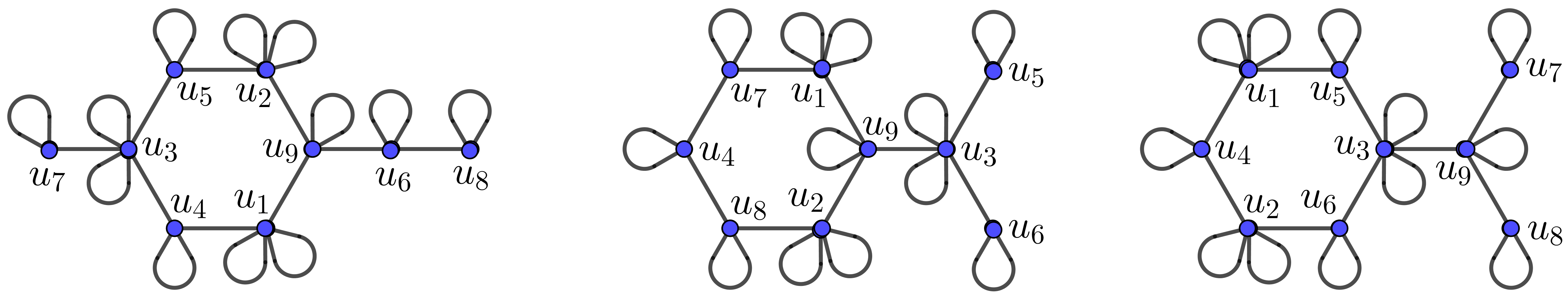} \\[6pt]
$c)$ \\[6pt]
\includegraphics[scale=0.65]{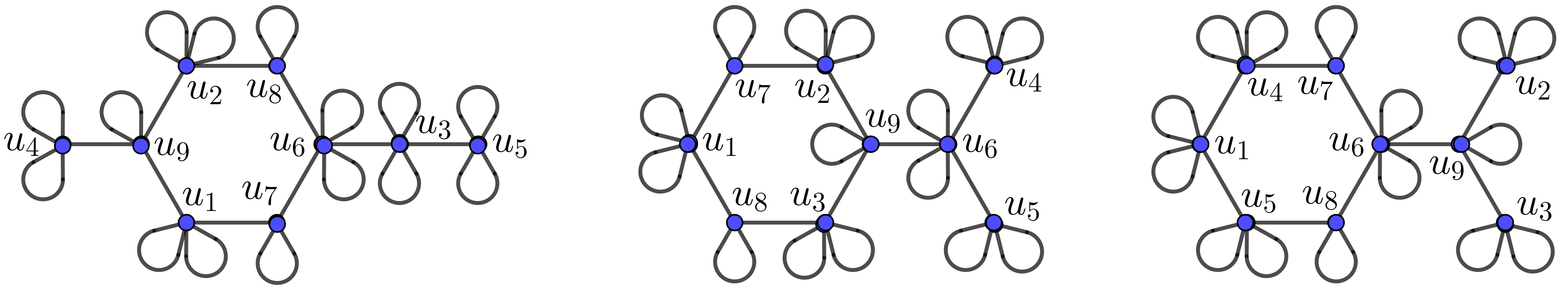} \\[6pt]
$d)$ 
{\small
\caption{
Four triplets of CMCs with 9-vertex unicyclic graphs as underlying graphs and:
a) the signature $(1,1,1)$,
b) the signature $(1,1,1,1,1,1)$,
c) the signature $(2,2,2,1,1,1,1,1,1)$, and
d) the signature $(2,2,2,2,2,2,1,1,1)$.
\label{fig-example-4}
}}
\end{figure}

\begin{figure}
\centering
\includegraphics[scale=0.65]{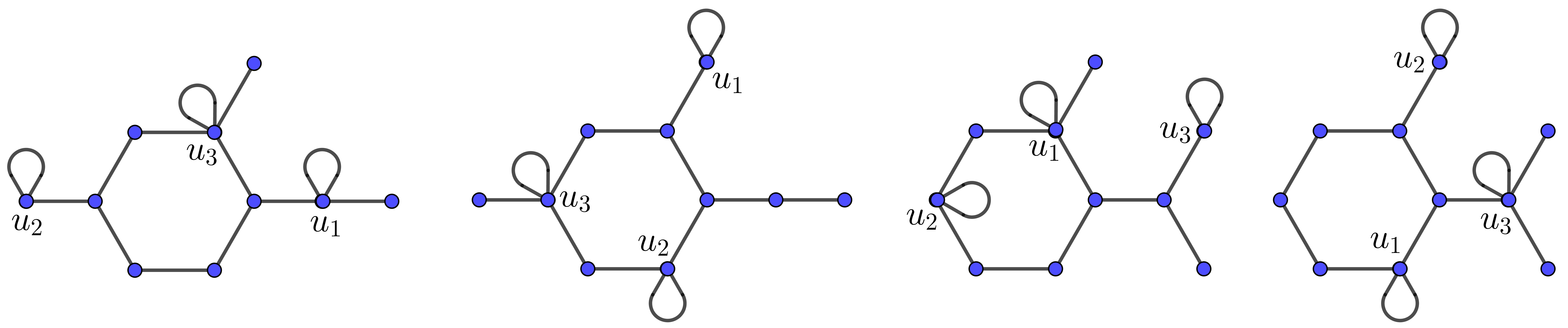} \\[6pt]
$a)$ \\[6pt]
\includegraphics[scale=0.65]{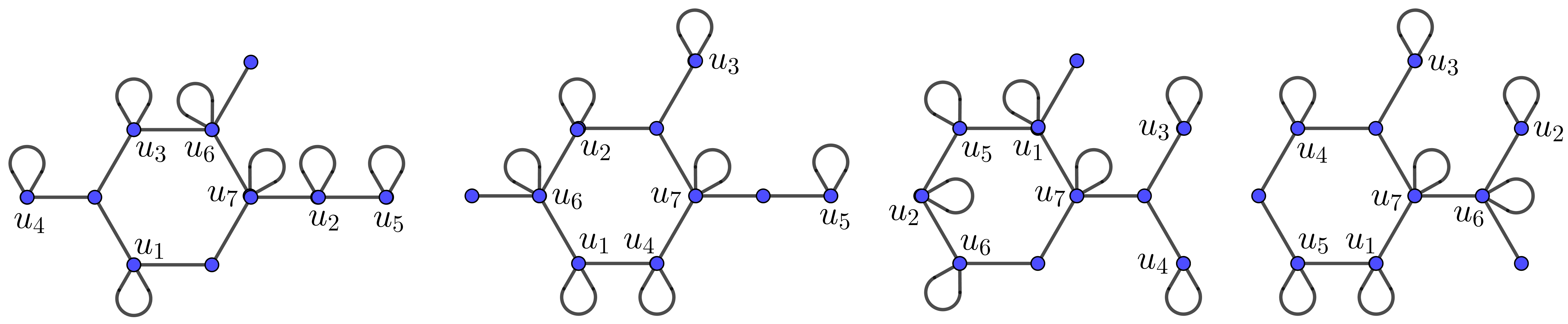} \\[6pt]
$b)$ 
{\small
\caption{Two quadruplets of CMCs 
with 10-vertex unicyclic graphs as underlying graphs and:
a) the signature $(1,1,1)$, 
b) the signature $(1,1,1,1,1,1,1)$.
\label{fig-example-8}
}}
\end{figure}

\begin{figure}
\centering
\includegraphics[scale=0.65]{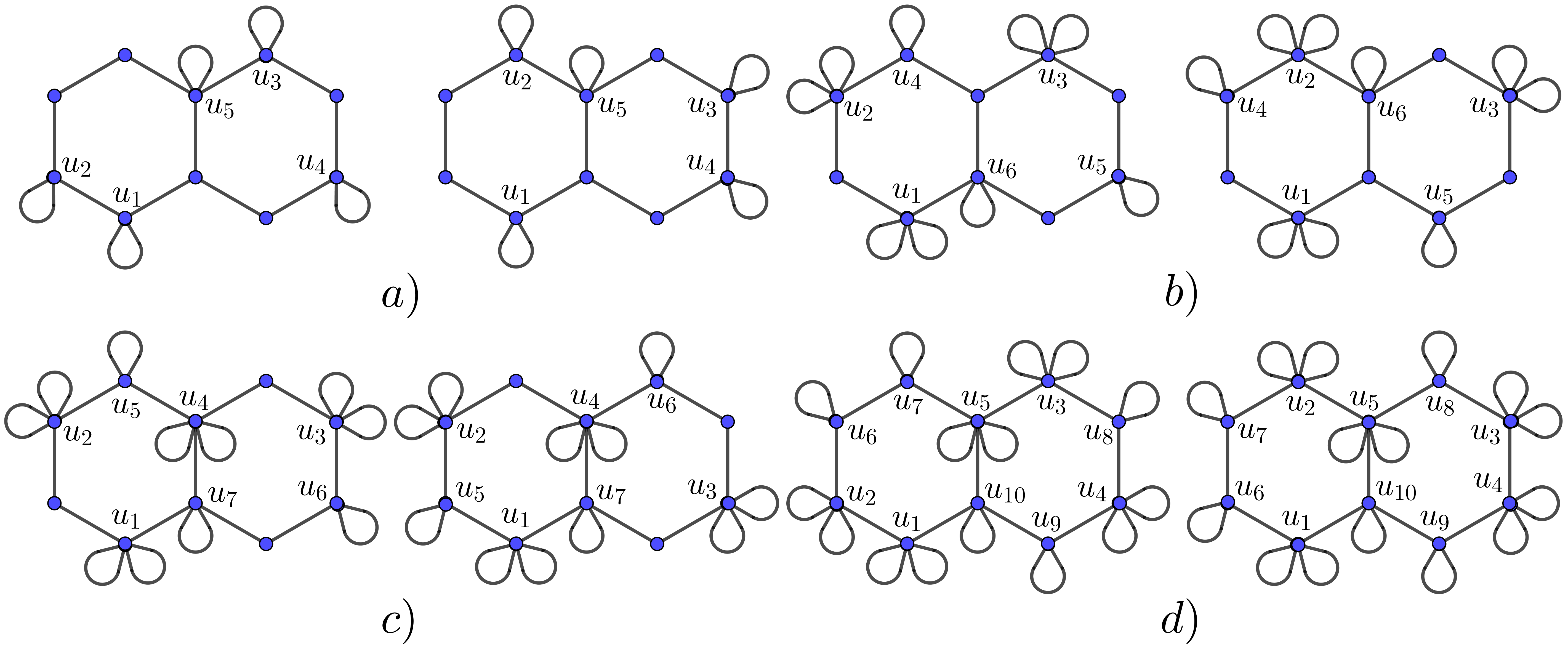}
{\small
\caption{Four pairs of CMCs 
with a benzenoid with two hexagons as an underlying graph and: 
a) the signature $(1, 1, 1, 1, 1)$,
b) the signature $(2, 2, 2, 1, 1, 1)$,
c) the signature $(2, 2, 2, 2, 1, 1, 1)$,
and d) the signature $(2, 2, 2, 2, 2, 1, 1, 1, 1, 1)$.
\label{fig-example-10}
}}
\end{figure}

\begin{figure}
\centering
\includegraphics[scale=0.9]{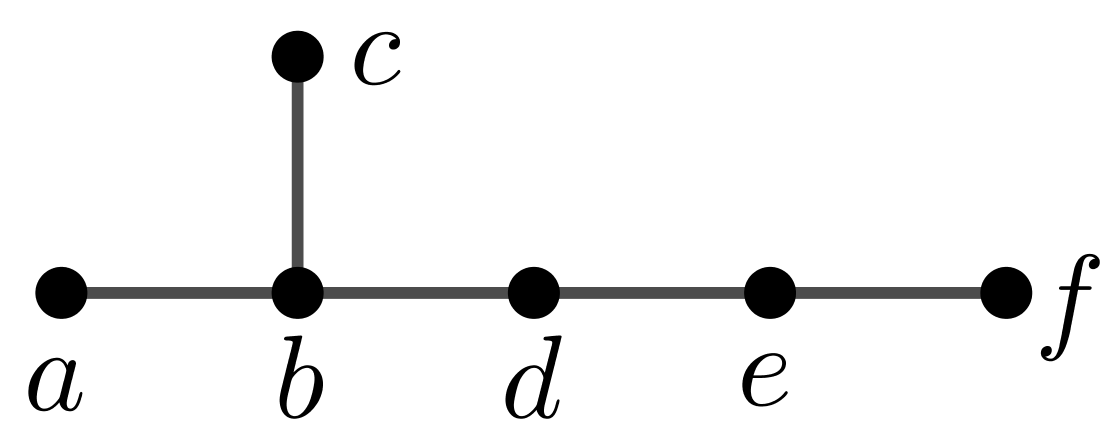} \\[12pt]
\begin{tabular}{lll}
\toprule
Signature & First vertex selection & Second vertex selection \\
\midrule
$(2, 1)$ & $(f, d)$ & $(a, e)$ \\
$(2, 1)$ & $(d, f)$ & $(a, b)$ \\
$(1, 1, 1)$ & $(d, a, e)$ & $(f, a, b)$ \\
$(2, 1, 1)$ & $(e, d, b)$ & $(b, f, e)$ \\
$(2, 1, 1)$ & $(d, e, b)$ & $(b, d, f)$ \\
$(2, 2, 2, 1, 1)$ & $(d, a, c, f, e)$ & $(f, a, c, d, b)$ \\
$(2, 2, 2, 1, 1)$ & $(f, a, c, e, b)$ & $(a, c, e, d, f)$ \\
$(2, 2, 2, 2, 1)$ & $(d, f, a, b, e)$ & $(a, c, e, b, d)$ \\
$(2, 2, 2, 2, 1)$ & $(d, f, a, e, b)$ & $(a, c, e, b, f)$ \\
$(2, 2, 2, 1, 1, 1)$ & $(d, a, e, f, c, b)$ & $(f, a, b, d, c, e)$ \\
\bottomrule
\end{tabular}
{\small
\caption{All ten pairs of CMCs with the tree on the top as an underlying graph
and the signatures with the maximum entry two.
Apart from the first, second and the last pair,
the remaining pairs were found by Hosoya in~\cite{hosoya-2}.
Further, there are 20 such pairs with the MSE~3,
30 pairs with MSE~4,
50 pairs with MSE~5, 
60 pairs with MSE~6,
90 pairs with MSE~7, and
110 pairs with MSE~8.
\label{fig-example-11}
}}
\end{figure}

We have run this search in a number of graph sets with various values of MSE.
Hosoya~\cite{hosoya-2} asked specifically for the new examples of CMCs 
in which the underlying graphs contain cycles or in which the underlying graphs are paths,
and such graph sets represent the majority of our search efforts.
The graph sets, MSEe and the numbers of examples of CMCs found
are shown in Table~\ref{tb-results}.
Detailed data about the examples found is available at~\cite{zenodo}.
The search in the sets of graphs with at most ten vertices was exhaustive.
The search among paths with between 11 and 20 vertices was incomplete,
but it still yielded sufficiently many examples to identify an infinite family of examples,
which will be described later.
Most of the examples consisted of the {\em pairs} of (underlying graph, signature, vertex subset) triplets
that satisfy the conditions~(\ref{eq-Hosoya-condition}).
However, 
the search has also found examples of triplets and even quadruplets of CMCs
with the underlying graphs among the connected graphs on seven vertices and the unicyclic graphs of girth six on 9 and 10 vertices.
These triplets and quadruplets, as well as four pairs with a benzenoid with two hexagons as an underlying graph,
are shown in Figs.~\ref{fig-example-1}--\ref{fig-example-10}.
As we can see from these figures,
the same underlying graphs often appear in several of these examples, 
but with different signatures and vertex selections.
Hosoya has already made this evident with seven examples of CMCs with MSE two
that all use the same 6-vertex tree shown in Fig.~\ref{fig-example-11} as an underlying graph.
Our exhaustive search yielded three new pairs of CMCs with MSE two
that are also based on this tree (see Fig.~\ref{fig-example-11}).
All this suggests that it is very likely that certain graphs will appear as underlying graphs
in a large (possibly infinite) number of CMCs.

\begin{table}
\centering
\begin{tabular}{p{7cm}cccc}
\toprule
Set of graphs & MSE & \# Pairs & \# Triplets & \# Quadruplets \\
                    &        &  of CMC  &  of CMC     &  of CMC \\
\midrule
Connected graphs, 5 vertices & 4 & 20 & & \\
Connected graphs, 6 vertices & 4 & 277 & & \\
Connected graphs, 7 vertices & 2 & 1215 & 3 & \\
Path, 8 vertices & 10 & 2788 & & \\
Path, 11 vertices & 1 & 4 & & \\
Path, 14 vertices & 1 & 10 & & \\
Path, 15 vertices & 1 & 3 & & \\
Path, 17 vertices & 1 & 11 & & \\
Path, 19 vertices & 1 & 10 & & \\
Path, 20 vertices & 1 & 9 & & \\
Trees with perfect matchings, 6 vertices & 3 & 14 & & \\
Trees with perfect matchings, 8 vertices & 2 & 89 & & \\
Trees with perfect matchings, 10 vertices & 1 & 105 & & \\
Unicyclic graphs, girth 6, 6 vertices & 2 & 1 & & \\
Unicyclic graphs, girth 6, 7 vertices & 2 & 2 & & \\
Unicyclic graphs, girth 6, 8 vertices & 2 & 52 & & \\
Unicyclic graphs, girth 6, 9 vertices & 2 & 745 & 4 & \\
Unicyclic graphs, girth 6, 10 vertices & 1 & 429 & 6 & 2 \\
Benzenoid, 2 hexagons & 2 & 4 & & \\
\bottomrule
\end{tabular}
\caption{Numbers of examples of cospectral multiple coalescences (CMC) found 
in selected sets of graphs for given maximum signature entry (MSE).}
\label{tb-results}
\end{table}

We paid special attention to finding examples of CMCs
in which the underlying graphs are paths.
Hosoya found one such example with the path~$P_8$ as an underlying graph and the signature $(1,1,1,1)$.
He posed the expectation that further such examples exist, 
with longer paths as underlying graphs and with larger signatures,
although he deemed their existence as rare~\cite[Section~5]{hosoya-2}.
In Table~\ref{tb-paths} we list the pairs of CMCs we found,
whose underlying graphs are paths with between 8 and~20 vertices and the MSE~1.
The exhaustive search with the path~$P_8$ as an underlying graph revealed
that there are 2788 pairs of CMCs with the MSE~10.
The searches with the paths $P_9$ and~$P_{10}$ as underlying graphs revealed no examples of CMCs.
The searches with the paths with between 11 and 20 vertices as underlying graphs and the MSE~1 were incomplete,
but they still uncovered a number of pairs, shown in Table~\ref{tb-paths}, 
that were sufficient to identify a rather general infinite family of examples of such CMCs.

\begin{table}
\centering
{\footnotesize
\begin{tabular}{rcll}
\toprule
Pair & Underlying path & First vertex selection & Second vertex selection \\
\midrule
1 & $P_8$ & $(0, 2, 3, 6)$ & $(0, 3, 5, 6)$ \\
\midrule
2 & $P_{11}$ & $(0, 3, 4, 8)$ & $(0, 4, 7, 8)$ \\
3 &               & $(0, 2, 3, 6, 9)$ & $(0, 3, 6, 8, 9)$ \\
4 &               & $(0, 2, 3, 5, 6, 9)$ & $(0, 3, 5, 6, 8, 9)$ \\
5 &               & $(0, 1, 3, 4, 5, 8, 9)$ & $(0, 1, 4, 5, 7, 8, 9)$ \\
\midrule
6 & $P_{14}$ & $(0, 4, 5, 10)$ & $(0, 5, 9, 10)$ \\
7 &               & $(1, 4, 6, 11)$ & $(1, 6, 9, 11)$ \\
8 &               & $(0, 2, 3, 6, 9, 12)$ & $(0, 3, 6, 9, 11, 12)$ \\
9 &               & $(0, 3, 5, 6, 9, 12)$ & $(0, 3, 6, 8, 9, 12)$ \\
10&               & $(0, 2, 4, 5, 7, 10, 12)$ & $(0, 2, 5, 7, 9, 10, 12)$ \\
11&               & $(0, 2, 3, 6, 8, 9, 12)$ & $(0, 3, 5, 6, 9, 11, 12)$ \\
12&               & $(0, 1, 4, 5, 6, 10, 11)$ & $(0, 1, 5, 6, 9, 10, 11)$ \\
13&               & $(0, 2, 3, 5, 6, 9, 12)$ & $(0, 3, 6, 8, 9, 11, 12)$ \\
14&               & $(0, 2, 3, 5, 6, 8, 9, 12)$ & $(0, 3, 5, 6, 8, 9, 11, 12)$ \\
15&               & $(0, 2, 3, 5, 6, 9, 11, 12)$ & $(0, 2, 3, 6, 8, 9, 11, 12)$ \\
\midrule
16&$P_{15}$ & $(0, 3, 4, 8, 12)$ & $(0, 4, 8, 11, 12)$ \\
17&               & $(0, 3, 4, 7, 8, 12)$ & $(0, 4, 7, 8, 11, 12)$ \\
18&               & $(0, 1, 3, 4, 5, 8, 9, 12, 13)$ & $(0, 1, 4, 5, 8, 9, 11, 12, 13)$ \\
\midrule
19&$P_{17}$ & $(0, 5, 6, 12)$ & $(0, 6, 11, 12)$ \\
20&               & $(1, 5, 7, 13)$ & $(1, 7, 11, 13)$ \\
21&               & $(1, 2, 5, 7, 8, 13, 14)$ & $(1, 2, 7, 8, 11, 13, 14)$ \\
22&               & $(0, 2, 5, 6, 8, 12, 14)$ & $(0, 2, 6, 8, 11, 12, 14)$ \\
23&               & $(0, 2, 3, 6, 9, 12, 15)$ & $(0, 3, 6, 9, 12, 14, 15)$ \\
24&               & $(0, 1, 5, 6, 7, 12, 13)$ & $(0, 1, 6, 7, 11, 12, 13)$ \\
25&               & $(0, 3, 5, 6, 9, 12, 15)$ & $(0, 3, 6, 9, 11, 12, 15)$ \\
26&               & $(0, 2, 3, 5, 6, 9, 12, 15)$ & $(0, 3, 6, 9, 11, 12, 14, 15)$ \\
27&               & $(0, 3, 5, 6, 8, 9, 12, 15)$ & $(0, 3, 6, 8, 9, 11, 12, 15)$ \\
28&               & $(0, 2, 3, 6, 9, 11, 12, 15)$ & $(0, 3, 5, 6, 9, 12, 14, 15)$ \\
29&               & $(0, 2, 3, 6, 8, 9, 12, 15)$ & $(0, 3, 6, 8, 9, 12, 14, 15)$ \\
\midrule
30&$P_{19}$ & $(1, 4, 6, 11, 16)$ & $(1, 6, 11, 14, 16)$ \\
31&               & $(0, 4, 5, 10, 15)$ & $(0, 5, 10, 14, 15)$ \\
32&               & $(0, 4, 7, 8, 12, 16)$ & $(0, 4, 8, 11, 12, 16)$ \\
33&               & $(0, 3, 4, 8, 12, 16)$ & $(0, 4, 8, 12, 15, 16)$ \\
34&               & $(1, 4, 6, 9, 11, 16)$ & $(1, 6, 9, 11, 14, 16)$ \\
35&               & $(0, 4, 5, 9, 10, 15)$ & $(0, 5, 9, 10, 14, 15)$ \\
36&               & $(0, 3, 4, 7, 8, 12, 16)$ & $(0, 4, 8, 11, 12, 15, 16)$ \\
37&               & $(0, 4, 7, 8, 12, 15, 16)$ & $(0, 3, 4, 8, 11, 12, 16)$ \\
38&               & $(0, 3, 4, 7, 8, 12, 15, 16)$ & $(0, 3, 4, 8, 11, 12, 15, 16)$ \\
39&               & $(0, 3, 4, 7, 8, 11, 12, 16)$ & $(0, 4, 7, 8, 11, 12, 15, 16)$ \\
\midrule
40&$P_{20}$ & $(0, 6, 7, 14)$ & $(0, 7, 13, 14)$ \\
41&               & $(1, 6, 8, 15)$ & $(1, 8, 13, 15)$ \\
42&               & $(2, 6, 9, 16)$ & $(2, 9, 13, 16)$ \\
43&               & $(1, 3, 6, 8, 10, 15, 17)$ & $(1, 3, 8, 10, 13, 15, 17)$ \\
44&               & $(1, 2, 6, 8, 9, 15, 16)$ & $(1, 2, 8, 9, 13, 15, 16)$ \\
45&               & $(0, 2, 6, 7, 9, 14, 16)$ & $(0, 2, 7, 9, 13, 14, 16)$ \\
46&               & $(0, 3, 6, 7, 10, 14, 17)$ & $(0, 3, 7, 10, 13, 14, 17)$ \\
47&               & $(0, 4, 6, 7, 11, 14, 18)$ & $(0, 4, 7, 11, 13, 14, 18)$ \\
48&               & $(0, 1, 6, 7, 8, 14, 15)$ & $(0, 1, 7, 8, 13, 14, 15)$ \\
\bottomrule
\end{tabular}
\caption{The pairs of CMCs with paths on 8 to 20 vertices as underlying graphs and the MSE~1
(hence the signature is $(1,\dots,1)$ in all these pairs).}
\label{tb-paths}
}
\end{table}

Several patterns could be initially observed from the pairs shown in Table~\ref{tb-paths}:
\begin{itemize}
\item[a)] In the pairs 1, 2, 6, 19 and 40, the number of vertices of the underlying path is $n=3k-1$,
              the first vertex selection is $\{0, k, 2k\}\cup\{k-1\}$ and the second vertex selection is $\{0, k, 2k\}\cup\{2k-1\}$.
\item[b)] In the pairs 7, 20 and 41, the number of vertices is $n=3k-1$,
              the first vertex selection is $\{1, k+1, 2k+1\}\cup\{k-1\}$ and 
              the second vertex selection is $\{1, k+1, 2k+1\}\cup\{2k-1\}$.
              Moreover, in the pair 42 
              the first vertex selection is $\{2, k+2, 2k+2\}\cup\{k-1\}$ and
              the second vertex selection is $\{2, k+2, 2k+2\}\cup\{k-1\}$.
\item[c)] In the pairs 1, 3, 8 and 23, the number of vertices is again $n=3k-1$,
              but the signature is of length~$k+1$ this time:
              the first vertex selection is $\{0,3,\dots,3k-3\}\cup\{2\}$ and
              the second vertex selection is $\{0,3,\dots,3k-3\}\cup\{3k-4\}$.
\item[d)] In the pairs 2, 16 and 33, the number of vertices is $n=4k-1$,
              the first vertex selection is $\{0,4,\dots,4k-4\}\cup\{3\}$ and
              the second vertex selection is $\{0,4,\dots,4k-4\}\cup\{4k-5\}$.
\end{itemize}

A slightly deeper analysis of these patterns reveals that they are actually all instances of a more general pattern,
in which the common part of two vertex selections consists of an initial sequence and its several translations,
with the remaining vertices from the vertex selections appearing at appropriate places between these translations.
In particular, we can prove the following cospectrality result of multiple coalescences 
with paths as underlying graphs and arbitrarily long signatures of ones.

\begin{theorem}
\label{th-underlying-paths}
For arbitrary integers $k\geq 3$, $m<k/2$, $d\geq 2$, and 
the integer tuple $(a_1,\dots,a_p)$ such that $0\leq a_1<\cdots<a_p\leq d-2$,
let $n=kd-1$, $v=md-1$, $w=(k-m)d-1$, and
let $u\colon\{1,\dots,k\}\times\{1,\dots,p\}\to\mathbb{Z}$ be defined as $u_{(i,j)}=(i-1)d+a_j$.
Then for any rooted graph~$G$ with the root~$r$
the multiple coalescences
\begin{equation}
\label{eq-path-v}
P_n(v=r)G(u_{(1,1)}=r)G(u_{(1,2)}=r)G\cdots(u_{(k,p)}=r)G
\end{equation}
and
\begin{equation}
\label{eq-path-w}
P_n(w=r)G(u_{(1,1)}=r)G(u_{(1,2)}=r)G\cdots(u_{(k,p)}=r)G
\end{equation}
are cospectral.
\end{theorem}

Before we proceed with the proof,
the following diagram serves to better illustrate the construction from the above theorem:
\begin{center}
\includegraphics[scale=0.75]{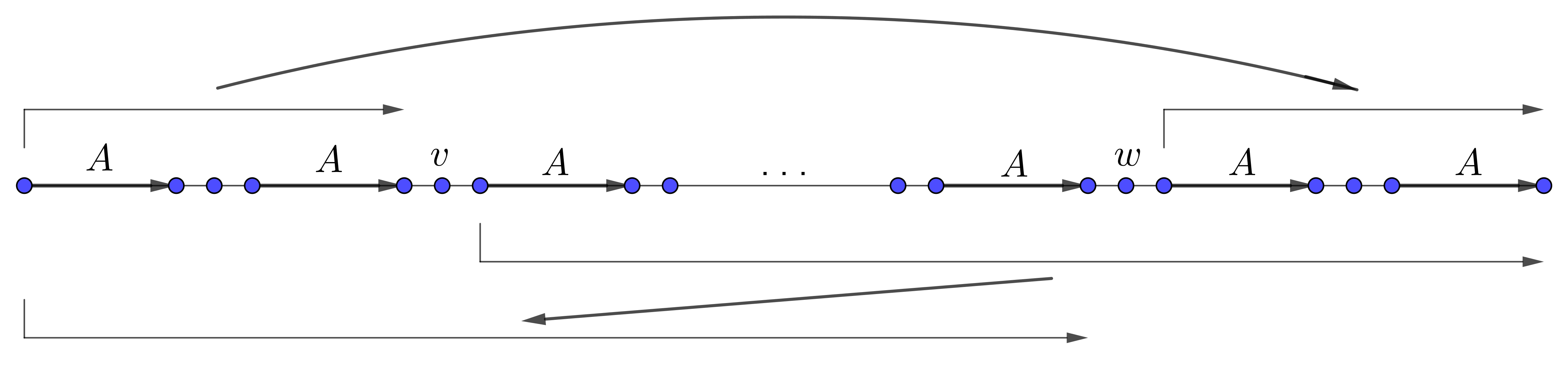}
\end{center}
The whole diagram represents the path~$P_n$,
where each of $k$ $A$-segments is a copy of the tuple $(a_1,\dots,a_p)$
translated for some multiple of~$d$ along the path.
An important consequence of the symmetric placement of vertices $v$ and~$w$ between the $A$-segments is that
$P_n-v$ and $P_n-w$ consist of two subpaths each
so that one has apparent bijections between the $A$-segments in these subpaths,
as illustrated above.
These bijections will enable us to easily prove that the conditions~(\ref{eq-Hosoya-condition})
hold for the vertex selections $(v,u_{(1,1)},\dots,u_{(k,p)})$ and $(w,u_{(1,1)},\dots,u_{(k,p)})$.

\bigskip\noindent
{\em Proof of Theorem~\ref{th-underlying-paths}.}\ \ 
The signature in both multiple coalescences (\ref{eq-path-v}) and~(\ref{eq-path-w}) 
consists of $kp+1$ ones.
From Corollary~\ref{co-Hosoya}, 
these multiple coalescences will be cospectral
if and only if for each $l=0,\dots,kp+1$ we have
\begin{align*}
  &\sum_{\substack{I\subseteq\{1,\dots,k\}\times\{1,\dots,p\}\\|I|=l-1}} P(P_n-v-\sum_{(i,j)\in I} u_{(i,j)}) 
  + \sum_{\substack{J\subseteq\{1,\dots,k\}\times\{1,\dots,p\}\\|J|=l}} P(P_n-\sum_{(i,j)\in J} u_{(i,j)}) \\
=&\sum_{\substack{I\subseteq\{1,\dots,k\}\times\{1,\dots,p\}\\|I|=l-1}} P(P_n-w-\sum_{(i,j)\in I} u_{(i,j)}) 
  + \sum_{\substack{J\subseteq\{1,\dots,k\}\times\{1,\dots,p\}\\|J|=l}} P(P_n-\sum_{(i,j)\in J} u_{(i,j)}),
\end{align*}
where in the sumations above we have separated the vertex subsets that contain $v$ or~$w$ 
from those that do not contain them.
The second sumations on both the left-hand side and the right-hand side are identical,
so it remains to prove that
\begin{equation}
\label{eq-path-subtract-v-and-w}
  \sum_{\substack{I\subseteq\{1,\dots,k\}\times\{1,\dots,p\}\\|I|=l}} P(P_n-v-\sum_{(i,j)\in I} u_{(i,j)}) 
=\sum_{\substack{I\subseteq\{1,\dots,k\}\times\{1,\dots,p\}\\|I|=l}} P(P_n-w-\sum_{(i,j)\in I} u_{(i,j)}) 
\end{equation}
for $l=0,\dots,kp$.
To show this, we construct the following bijection~$\theta$ 
on the set $\{0,\dots,kd-2\}\setminus\{v,w\}$ of vertices of $P_n-v-w$,
which formalizes the bijection illustrated in the diagram above:
$$
\theta x = \left\{\begin{array}{rl}
x + (k-m)d, & \mbox{if }x<v, \\
x - md, & \mbox{if }x>v.
\end{array}\right.
$$
Now, the graph $P_n-v$ consists of 
the subpath $V_l$ on vertices labeled $0,\dots,md-2$ of length $md-2$ and
the subpath $V_r$ on vertices labeled $md,\dots,kd-1$ of length $(k-m)d-2$,
while the graph $P_n-w$ consists of
the subpath $W_r$ on vertices labeled $(k-m)d,\dots,kd-1$ of length $md-2$ and
the subpath $W_l$ on vertices labeled $0,\dots,(k-m)d-2$ of length $(k-m)d-2$.
It is straightforward to see that
$\theta|_{V_l}$ is an isomorphism from $V_l$ to $W_r$,
$\theta|_{V_r}$ is an isomorphism from $V_r$ to $W_l$,
and that
$$
\theta u_{(i,j)} = \left\{\begin{array}{rl}
u_{(i+k-m,j)} & \mbox{if }i\leq m, \\
u_{(i-m,j)}, & \mbox{if }m<i.
\end{array}\right.
$$
Let us further define $\theta'$ on $\{1,\dots,k\}\times\{1,\dots,p\}$ by
$$
\theta'(i,j) = \left\{\begin{array}{rl}
(i+k-m,j) & \mbox{if }i\leq m, \\
(i-m,j), & \mbox{if }m<i,
\end{array}\right.
$$
so that $\theta u_{(i,j)}=u_{\theta'(i,j)}$.
Then 
\begin{align*}
P_n-v-\sum_{(i,j)\in I} u_{(i,j)} 
& \cong P_n-w-\sum_{(i,j)\in I} \theta u_{(i,j)} \\
& = P_n-w-\sum_{(i,j)\in I} u_{\theta'(i,j)} \\
& = P_n-w-\sum_{(i',j')\in \theta'I} u_{(i',j')},
\end{align*}
and consequently
$$
P(P_n-v-\sum_{(i,j)\in I} u_{(i,j)}) = P(P_n-w-\sum_{(i',j')\in \theta'I} u_{(i',j')}).
$$
Since $\theta'$ is a permutation of $\{1,\dots,k\}\times\{1,\dots,p\}$,
we have that
$\theta'I$ ranges through all the $l$-element subsets of $\{1,\dots,k\}\times\{1,\dots,p\}$
when $I$ ranges through all such subsets for any fixed $l$, $0\leq l\leq kp$.
Hence 
\begin{align*}
    \sum_{\substack{I\subseteq\{1,\dots,k\}\times\{1,\dots,p\}\\|I|=l}} P(P_n-v-\sum_{(i,j)\in I} u_{(i,j)})
=&\sum_{\substack{I\subseteq\{1,\dots,k\}\times\{1,\dots,p\}\\|I|=l}} P(P_n-w-\sum_{(i',j')\in \theta'I} u_{(i',j')}) \\
=&\sum_{\substack{J\subseteq\{1,\dots,k\}\times\{1,\dots,p\}\\|J|=l}} P(P_n-w-\sum_{(i',j')\in J} u_{(i',j')}),
\end{align*}
where $J$ denotes $\theta'I$ in the last equality.
This proves~(\ref{eq-path-subtract-v-and-w}),
and consequently shows that the multiple coalescences (\ref{eq-path-v}) and~(\ref{eq-path-w}) are cospectral.
\hfill$\Box$

One can easily inspect that a large percentage of the parameters of CMC pairs collected in Table~\ref{tb-paths}
are particular instances of the infinite family of such parameters identified in Theorem~\ref{th-underlying-paths}.
The parameters of the remaining CMC pairs from Table~\ref{tb-paths},
and especially the large number of pairs of CMCs with $P_8$ as the underlying graph and the MSE 10,
suggest that it would be quite possible to identify further infinite families of CMCs pairs with paths as underlying graphs.

\section*{Acknowledgements}

This work was supported and funded by Kuwait University Research Grant No. SM03/20.

\end{document}